\documentclass[a4paper]{article}
\usepackage[utf8]{inputenc}

\usepackage[english]{babel}
\usepackage{graphicx}
\usepackage{xcolor}
\usepackage{amsmath,amsthm}
\usepackage{tikz}
\usetikzlibrary{tikzmark}
\usepackage{amssymb}
\usepackage{mathrsfs}
\usepackage{mathabx}
\usepackage{txfonts,pxfonts,tikz}
\usepackage{tgschola}
\usepackage[T1]{fontenc}
\usepackage[margin= 1 in]{geometry}
\usepackage{mathtools}
\usepackage{enumitem}   
\usepackage{listings}

\newcommand{\defeq}{\vcentcolon=}

\usepackage[backend=biber,style=nature]{biblatex}
\usepackage{hyperref}
\hypersetup{
    colorlinks=true,
    citecolor=cyan,
    }

\newcommand{\Z}{\mathbb{Z}}

\title{Existential Inclusions of Bi-exact Groups are Conjugacy Representation Rigid}
\author{Connor MacMahon}

\theoremstyle{definition}

\theoremstyle{plain}

\begin{document}

\maketitle

\begin{abstract}
    If $\Lambda$ is a non-amenable bi-exact group and $\Lambda \hookrightarrow \Gamma$ is an existential embedding, then each of the intersections $\Lambda \cap g \Lambda g^{-1}$ for $g$ a member of $\Gamma \backslash \Lambda$ is amenable. This in conjunction with work of Bekka and Kalantar \cite{BK} demonstrates that in this situation, the weak equivalence class of the quasi-regular representation $\lambda_{\Gamma/\Lambda}$ determines $\Lambda$ up to conjugacy among the self-commensurating subgroups of $\Gamma$. 
\end{abstract}

\section{Introduction}

Suppose that $\Gamma$ is a countable discrete group, and that $\Lambda$ is a subgroup which belongs to a specified class $\mathscr{C}$ of subgroups of $\Gamma$. Originating in the area of spectral geometry, it has become a question of interest to compare the quasi-regular representation $\lambda_{\Gamma/\Lambda}$ with those arising from other subgroups in the class $\mathscr{C}$. 

\medskip

Indeed, Sunada \cite{Sunada} leveraged the comparison of quasi-regular representations in the setting of $\Gamma$ a finite group in order to produce two isospectral Riemannian manifolds which are not isometric.

\medskip

Recently, work in this situation was revived by Bekka and Kalantar \cite{BK}. In this work, they studied subgroups $\Lambda$ of a countable discrete group $\Gamma$ which exhibit a rigidity phenomenon involving the quasi-regular representation. In particular, $\Lambda$ is conjugacy representation rigid in the class $\mathscr{C}$ provided that if $\Lambda'$ is another group belonging to $\mathscr{C}$ and $\lambda_{\Gamma/\Lambda}$ is weakly equivalent to $\lambda_{\Gamma/\Lambda'}$, then $\Lambda$ and $\Lambda'$ are conjugate within $\Gamma$ (I am referring to their rigidity notion as "conjugacy representation rigidity" in order to avoid ambiguity regarding other notions of representation rigidity in the literature).

\medskip

Their results concern the situation in which the class $\mathscr{C}$ consists of self-commensurating subgroups of $\Gamma$. This is a particularly relevant setup, because the quasi-regular representation associated to a self-commensurating subgroup is irreducible. 

\medskip

Suppose that $\Lambda$ is a non-amenable subgroup of $\Gamma$ such that $\Lambda \cap g\Lambda g^{-1}$ is amenable for each member $g$ of $\Gamma \backslash \Lambda$. In this situation, $\Lambda$ is referred to as a-normal \cite{BK} (a-normality does not formally require non-amenability of $\Lambda$, but the amenable case is uninteresting). Then $\Lambda$ is self-commensurating, and Bekka and Kalantar proved that $\Lambda$ is conjugacy representation rigid in the class of self-commensurating subgroups of $\Gamma$ \cite{BK}.

\medskip

The purpose of this brief note is to provide further examples of this situation arising from first order logic. The main mechanism for this is the following Theorem.

\medskip

\noindent \textbf{Theorem:} Let $\Lambda$ denote a bi-exact group and suppose that $\Lambda \hookrightarrow \Gamma$ is an existential embedding. Then $\Lambda$ is an a-normal subgroup of $\Gamma$.

\medskip

In work of Perin \cite{FO}, there are given many examples of elementary inclusions $\Lambda \hookrightarrow \Gamma$ with $\Gamma$ torsion free hyperbolic in which $\Lambda$ is actually malnormal in $\Gamma$ (i.e. the intersections $\Lambda \cap g\Lambda g^{-1}$ are trivial for $g$ a member of $\Gamma \backslash \Lambda$). Recall that by the permanence properties of existential and universal formulas, if $\Lambda \leq \Gamma$ is an elementary embedding, and
\begin{align*}
    \Lambda \leq K \leq \Gamma,
\end{align*}
then the inclusion $\Lambda \leq K$ is existential. Moreover, malnormality of $\Lambda$ in $\Gamma$ evidently implies malnormality of $\Lambda$ within $K$. Therefore, Perin's work provides numerous examples of the above phenomenon for torsion free hyperbolic groups, such that the main Theorem is partially covered by existing work.

\medskip
\medskip

\noindent \textbf{Acknowledgement:} The author acknowledges support from NSF CAREER award DMS-21447.

\section{Preliminaries on First Order Logic}

Here, the necessary first order logic is introduced following Tent and Ziegler \cite{ModTHy}. Within the language of groups, the terms are words in the variables denoting group elements, their inverses, and the identity element. An atomic formula consists of two terms separated by "$=$".

\medskip

To describe first order formulae, one combines atomic formulae, logical connectors ($\vee$, $\wedge$, $\neg$) and quantifiers ($\forall, \exists$). Where no free variables are present (i.e. each variable is quantified over), a formula is referred to as a sentence. For example, the following is a sentence in the first order language of groups:
\begin{align*}
    \forall g_1 \exists g_2 (g_1g_2g_1^{-1}g_2^{-1}=1) \wedge \neg(g_2 =1).
\end{align*}
The \textit{elementary theory} of a group $G$ is the collection of first order sentences which hold true in $G$. An embedding $H \hookrightarrow G$ of groups is called elementary if each first order sentence $\phi$ in the language of groups with constants in $H$ satisfies $G \models \phi$ if and only if $H \models \phi$.

\medskip

Observe that if $H$ is elementarily embedded within $G$, then $H$ and $G$ have the same elementary theory. Indeed, a first order sentence $\phi$ in the language of groups is also a first order sentence in the language of groups with constants in $H$. If $G$ and $H$ have the same elementary theory, then they are called elementarily equivalent. More precisely, $G$ and $H$ are elementarily equivalent provided that for each first order sentence $\phi$ in the language of groups, $G \models \phi$ if and only if $H \models \phi$.

\medskip

It should be noted that $H$ being elementarily embedded within $G$ is strictly stronger than $H$ being elementarily equivalent to $G$ (for instance, $2\Z \subset \Z$ is not an elementary embedding. See Guirardel, Levitt and Sklinos \cite{FO} 3.1.1). 

\medskip

The main focus below is existential embeddings. An existential sentence is a first order sentence $\phi$ in the language of groups which takes the following form:
\begin{align*}
    \exists x_1 \exists x_2 \dots \exists x_n \psi(x_1,\dots,x_n),
\end{align*}
where $\psi$ is a quantifier free formula. The \textit{existential theory} of a group $G$ consists of the collection of existential sentences which hold true in $G$. As with elementary embeddings, an embedding $H \hookrightarrow G$ of groups is called existential provided each existential sentence $\phi$ in the language of groups with constants in $H$ satisfies $G \models \phi$ if and only if $H \models \phi$. 

\section{Preliminaries on Bi-exactness}
The notion of a bi-exact group was introduced by Ozawa in 2004 \cite{Kurosh} towards deformation/rigidity results in the theory of II$_1$ factors. 

\medskip

A group $\Lambda$ is bi-exact provided the $\Lambda \times \Lambda$ action by left/right translations on the Stone-Čech remainder $\partial^\beta \Lambda \defeq \beta \Lambda \backslash \Lambda$ is amenable. Amenability of the action $\Lambda \times \Lambda \curvearrowright \partial^\beta \Lambda$ asserts the existence of a sequence of continuous maps $\mu^n : \partial^\beta \Lambda \rightarrow \text{Prob}(\Lambda \times \Lambda)$ such that:
\begin{align*}
   \lim_{n \rightarrow \infty} \sup_{x \in \partial^\beta \Lambda} \|s \cdot \mu^n_{x} - \mu^n_{s \cdot x}\|_1 = 0
\end{align*}
for each $s$ belonging to $\Lambda \times \Lambda$ \cite{Kurosh}. Here, 
\begin{align*}
    \text{Prob}(\Lambda \times \Lambda) \defeq \{\mu \in \ell^1(\Lambda \times \Lambda) \, | \, \mu \geq 0 \, \text{and} \, \|\mu\|_1=1\}
\end{align*}
and $(s \cdot \mu)(t) = \mu(s^{-1}t)$. For more equivalent conditions defining bi-exactness, see Brown and Ozawa \cite{BOzawa}. The importance of this notion in the present context comes from the following result, which is present in various forms in the literature.

\medskip

\noindent \textbf{Proposition (\cite{Comment}, \cite{DingPete}):} Suppose $\Lambda$ is a bi-exact group. Then for each non-principal ultrafilter $\omega$ on $\mathbb{N}$ there holds:
\begin{align*}
    _{L(\Lambda)}[L^2(L(\Lambda)^\omega) \ominus \ell^2(\Lambda)]_{L(\Lambda)} \prec \ell^2(\Lambda) \otimes \ell^2(\Lambda).
\end{align*}
That is, the $L(\Lambda)$-$L(\Lambda)$ bimodule on the left hand side is weakly coarse.
\begin{proof}
  A proof by literature. It was demonstrated by Ozawa \cite{Comment} that the group von Neumann algebra of a bi-exact group satisfies $W^*$AO. Now, apply Theorem 7.20 of Ding and Peterson \cite{DingPete}.
\end{proof}

The class of bi-exact groups is sizeable, containing all word hyperbolic groups, discrete subgroups of connected simple Lie groups of rank $1$, small cancellation groups, amenable groups, and the semidirect product $\Z^2 \rtimes SL_2(\Z)$ (\cite{BOzawa}, \cite{Small}). Moreover, this class is stable under measure equivalence \cite{ME}, $W^*$-equivalence \cite{DingPete}, free products with finite amalgamation \cite{Kurosh}, and passing to subgroups \cite{Kurosh}.

\section{Preliminaries on Representation Rigidity}

Following the treatment of Bekka and Kalantar \cite{BK}, the setup for conjugacy representation rigidity is introduced. Suppose that $\Lambda$ is a subgroup of a group $\Gamma$. Define the commensurator subgroup $\text{Comm}_\Gamma(\Lambda)$ to consist of those $g$ in $\Gamma$ such that $\Lambda$ and $g\Lambda g^{-1}$ are commensurable. That is, those $g$ for which the intersection $\Lambda \cap g\Lambda g^{-1}$ has finite index in both $\Lambda$ and $g \Lambda g^{-1}$. 

\medskip

Observe that $\Gamma$ acts by left multiplication on the coset space $\Gamma/\Lambda$. This action induces a unitary representation $\lambda_{\Gamma/\Lambda}$ on $\ell^2(\Gamma/\Lambda)$ called the quasi-regular representation associated to $\Gamma$. By a Theorem of Mackey (see \cite{BK} 2.6), if $\Lambda$ is self-commensurating (i.e. $\Lambda = \text{Comm}_\Gamma(\Lambda)$), then $\lambda_{\Gamma/\Lambda}$ is an irreducible representation of $\Gamma$. 

\medskip

It is then a natural question to ask when two self-commensurating subgroups $\Lambda_1,\Lambda_2$ of $\Gamma$ define unitarily (resp. weakly) equivalent irreducible representations of $\Gamma$.

\medskip

Towards an answer to this question, Bekka and Kalantar \cite{BK} defined the class $\text{Sub}_{\text{sg}}(\Gamma)$ of subgroups of $\Gamma$ with spectral gap. This class consists of those subgroups $\Lambda$ of $\Gamma$ such that the trivial representation $1_\Lambda$ is isolated in the spectrum of the natural representation of $\Lambda$ on $\ell^2(\Gamma/\Lambda)$. They then demonstrate that self-commensurating subgroups in this class are conjugacy representation rigid in the following sense.

\medskip

\noindent \textbf{Theorem A (\cite{BK}):} Suppose $\Lambda_1$ is a self-commensurating subgroup of a countable group $\Gamma$. If $\Lambda_1$ has spectral gap, then for all self-commensurating subgroups $\Lambda_2$ of $\Gamma$, the quasi-regular representations $\lambda_{\Gamma/\Lambda_1}$ and $\lambda_{\Gamma/\Lambda_2}$ are weakly equivalent if and only if $\Lambda_1$ and $\Lambda_2$ are conjugate subgroups of $\Gamma$. 

\medskip

If $\Lambda_1$ is a self-commensurating subgroup of $\Gamma$ which satisfies the conclusion of this Theorem, it is said to be conjugacy representation rigid in the class of self-commensurating subgroups of $\Gamma$. In applying this Theorem, Bekka and Kalantar further introduce a sufficient condition on subgroups which guarantees the spectral gap property. It requires the following definition.

\medskip

\noindent \textbf{Definition:} Suppose $\Lambda$ is a subgroup of a countable group $\Gamma$. Then $\Lambda$ is a-normal provided $\Lambda \cap g\Lambda g^{-1}$ is amenable for each member $g$ of $\Gamma \backslash \Lambda$. 

\medskip

In their paper, Bekka and Kalantar (\cite{BK} 4.11) demonstrate that if a non-amenable subgroup $\Lambda$ of a countable group $\Gamma$ is a-normal, then it has the spectral gap property. Observe immediately that a non-amenable a-normal subgroup $\Lambda$ of $\Gamma$ is self-commensurating, such that in this case the spectral gap translates via the above Theorem into a representation rigidity result. In fact, such subgroups $\Lambda$ are \textit{strongly self-commensurating} in the sense that the intersections $\Lambda \cap g \Lambda g^{-1}$ for $g$ a member of $\Gamma \backslash \Lambda$ each have infinite index in $\Lambda$. 

\section{The Main Theorem}

Using a trick identical to that in work of Hayes, Jekel, and Kunnawalkam Elayavalli (\cite{RandEmbed} Proposition 3.1) along with some basic representation theory, the main theorem is proved.

\medskip

\noindent \textbf{Theorem:} Let $\Lambda$ denote a bi-exact group and suppose that $\Lambda \hookrightarrow \Gamma$ is an existential embedding. Then $\Lambda$ is an a-normal subgroup of $\Gamma$.
\begin{proof}
    Set $M=L(\Lambda)$, the group von Neumann algebra. By the bi-exactness of $\Lambda$ and the above Proposition,
    \begin{align*}
        L^2(M^\omega) \ominus \ell^2(\Lambda) \prec \ell^2(\Lambda) \otimes \ell^2(\Lambda)
    \end{align*}
    for each non-principal ultrafilter $\omega$ on $\mathbb{N}$. Now, suppose that $\Lambda \hookrightarrow \Gamma$ is an existential embedding. As was demonstrated in work of Kunnawalkam Elayavalli \cite{Sri}, there exists a non-principal ultrafilter $\omega$ on $\mathbb{N}$ and an embedding $L(\Gamma) \hookrightarrow M^\omega$ extending the diagonal embedding of $M$:
    \begin{align*}
        M \leq L(\Gamma) \leq M^\omega.
    \end{align*}
    Therefore,
    \begin{align*}
        {_M}[\ell^2(\Gamma)\ominus  \ell^2(\Lambda)]_M \prec L^2(M^\omega) \ominus \ell^2(\Lambda),
    \end{align*}
    which is weakly coarse by the above. Fix $g$ in $\Gamma \backslash \Lambda$. Complete $g \defeq g_1$ to a system $\{g_i\}_{i \in \mathbb{N}}$ of double coset representatives such that $\Gamma \backslash \Lambda$ is divided as the disjoint union:
    \begin{align*}
        \Gamma \backslash \Lambda = \bigcup_i \Lambda g_i \Lambda.
    \end{align*}
    Therefore,
    \begin{align*}
        {_M}[\ell^2(\Gamma)\ominus  \ell^2(\Lambda)]_M  \cong \bigoplus_{i} \ell^2(\Lambda g_i \Lambda)
    \end{align*}
    as $M-M$ bimodules. In particular, $\ell^2(\Lambda g\Lambda)$ is weakly coarse. Now, observe that $\Lambda \times \Lambda$ acts transitively on the double coset $\Lambda g\Lambda$. Therefore, if $H$ is the stabilizer of $g$, then:
    \begin{align*}
        \ell^2(\Lambda g \Lambda) \cong \ell^2(\Lambda \times \Lambda/H) = \lambda_{\Lambda \times \Lambda/H},
    \end{align*}
    the quasi-regular representation. Observe here that the $M$-$M$ bimodule structure translates into a unitary representation of $\Lambda \times \Lambda$ as the map $(x,y) \mapsto (x,y^{-1})$ is an isomorphism $\Lambda \times \Lambda^\text{op} \cong \Lambda \times \Lambda$. Under this identification, the coarse bimodule corresponds to the left regular representation. That $\ell^2(\Lambda g \Lambda)$ is weakly coarse means that $\lambda_{\Lambda \times \Lambda/H} \prec \lambda_{\Lambda \times \Lambda}$. Recall that $\lambda_{\Lambda \times \Lambda/H} = \text{Ind}_H^{\Lambda \times \Lambda} 1_H$, such that standard induction/restriction arguments (see Bekka and Kalantar \cite{BK} Corollary 2.3) yield $1_H \prec \lambda_H$. Since the trivial representation is weakly contained in the left regular representation, $H$ is amenable. Now, fix $(x,y)$ in $H$. By definition,
    \begin{align*}
        xgy^{-1} = g.
    \end{align*}
    Rearranging, $y = g^{-1}xg$. Hence, both $x$ and $g^{-1}x g$ belong to $\Lambda$, such that $x$ is a member of $\Lambda \cap g\Lambda g^{-1}$. Conversely, it is easily seen that a pair of the form $(x,g^{-1}xg)$ for $x$ a member of $\Lambda \cap g\Lambda g^{-1}$ belongs to $H$. In short, $H \cong \Lambda \cap g\Lambda g^{-1}$, so that $\Lambda \cap g\Lambda g^{-1}$ is amenable. Since $g$ was an arbitrary member of $\Gamma \backslash \Lambda$, the subgroup $\Lambda$ is a-normal.
\end{proof}

Using bi-exactness alone, the conclusion of a-normality is the strongest possible. Indeed, this is demonstrated by the case in which both $\Gamma$ and $\Lambda$ are abelian.

\medskip

\noindent \textbf{Remark:} Observe that the bimodule $\ell^2(\Lambda g \Lambda)$ being weakly coarse gave rise to the amenability of the intersection $\Lambda \cap g \Lambda g^{-1}$. It should be noted that if the bimodules $\ell^2(\Lambda g \Lambda)$ for $g$ a member of $\Gamma \backslash \Lambda$ are coarse (i.e. embed into a direct sum of the coarse bimodule $\ell^2(\Lambda) \otimes \ell^2(\Lambda)$), then $\Lambda$ is seen to be almost malnormal (\cite{Consequences} section 2). That is, the intersections $\Lambda \cap g \Lambda g^{-1}$ for $g$ a member of $\Gamma \backslash \Lambda$ are finite. 

\medskip

Now, the representation rigidity Corollary is immediate.

\medskip

\noindent \textbf{Corollary:} Let $\Lambda$ be a non-amenable bi-exact group and $\Lambda \hookrightarrow \Gamma$ an existential embedding. Both $\Lambda$ and $\Gamma$ are assumed to be countable and discrete. Then $\Lambda$ is strongly self-commensurating within $\Gamma$ and $\Lambda$ is conjugacy representation rigid in the class of self-commensurating subgroups of $\Gamma$.
\begin{proof}
    Apply Bekka and Kalantar \cite{BK} Theorem A and the main Theorem. 
\end{proof}

\medskip

While not directly applicable to representation rigidity, it should be noted that another Theorem of Bekka and Kalantar applies in this setting as well. This yields a Corollary of independent interest regarding the ideal structure in the $C^*$-algebra generated by the quasi-regular representation.

\medskip

\noindent \textbf{Corollary:} Let $\Lambda$ be a non-amenable bi-exact group, and $\Lambda \hookrightarrow \Gamma$ an existential embedding. Then $C^*_{\lambda_{\Gamma/\Lambda}}(\Gamma)$ contains the ideal $K$ of compact operators on $\ell^2(\Gamma/\Lambda)$, and $K$ is minimal in the sense that it is contained within each nonzero, closed, two-sided ideal of $C^*_{\lambda_{\Gamma/\Lambda}}(\Gamma)$. 
\begin{proof}
    Apply Bekka and Kalantar \cite{BK} Theorem B and the main Theorem.
\end{proof}

\printbibliography

\end{document}